\pdfoutput=1
\documentclass[twocolumn,10pt]{asme2ej}

\usepackage{amsmath}
\usepackage{amssymb}
\usepackage{subfigure}
\usepackage[normalem]{ulem}
\usepackage{graphicx}
\usepackage{epstopdf}
\usepackage{tikz}
\usetikzlibrary{arrows,shapes}

\newcommand{\E}{\mathbb{E}}

\newcommand{\samp}{\theta}
\newcommand{\loss}{\bar{\omega}}
\newcommand{\turn}{\Delta\beta}

\newcommand{\des}{\mathbf{d}}

\newcommand{\tol}{\boldsymbol\sigma}

\newcommand{\bxi}{\boldsymbol\xi}

\title{The Implications of Tolerance Optimization on Compressor Blade Design}

\author{Eric A. Dow
    \affiliation{
	Graduate Student\\
	Aerospace Computational Design Laboratory\\
	Department of Aeronautics and Astronautics\\
	Massachusetts Institute of Technology\\
	Cambridge, Massachusetts 02139\\
    Email: ericdow@mit.edu
    }	
}

\author{Qiqi Wang
    \affiliation{ 
    Assistant Professor\\
	Aerospace Computational Design Laboratory\\
	Department of Aeronautics and Astronautics\\
	Massachusetts Institute of Technology\\
	Cambridge, Massachusetts 02139\\
    Email: qiqi@mit.edu
    }
}

\begin{document}

\maketitle    

%%%%%%%%%%%%%%%%%%%%%%%%%%%%%%%%%%%%%%%%%%%%%%%%%%%%%%%%%%%%%%%%%%%%%%

\begin{abstract}
{\it Geometric variability increases performance variability and degrades the mean performance of turbomachinery compressor blades. These detrimental effects can be reduced by using robust optimization to design the blade geometry or by imposing stricter manufacturing tolerances. This paper presents a novel computational framework for optimizing compressor blade manufacturing tolerances, and incorporates this framework into existing robust geometry design frameworks. Optimizations of an exit guide vane geometry are conducted. The single-point optimal geometry is found to depend on the manufacturing tolerances due to a switch in the dominant loss mechanism. Multi-point geometry optimization avoids this switch so that the geometry and tolerance optimization problems are decoupled.

}
\end{abstract}

%%%%%%%%%%%%%%%%%%%%%%%%%%%%%%%%%%%%%%%%%%%%%%%%%%%%%%%%%%%%%%%%%%%%%%

\begin{nomenclature}
\entry{$e$}{manufacturing error field}
\entry{$s$}{surface location}
\entry{$\samp$}{element of sample space}
\entry{$C$}{covariance function}
\entry{$x$}{coordinates of blade surface}
\entry{$\hat{n}$}{surface normal}
\entry{$\rho$}{correlation function}
\entry{$\sigma$}{standard deviation}
\entry{$L$}{correlation length}
\entry{$\lambda$}{Karhunen-Lo\`{e}ve eigenvalue}
\entry{$\phi$}{Karhunen-Lo\`{e}ve eigenfunction}
\entry{$\xi$}{Karhunen-Lo\`{e}ve random variable coefficient}
\entry{$\loss$}{loss coefficient: $(p_{o2} - p_{o1})/(p_{o1} - p_1)$}
\entry{$\des$}{vector of nominal geometry design variables}
\entry{$\beta$}{flow angle}
\entry{$\alpha$}{incidence angle: $\beta - \beta_{des}$}
\entry{$C_p$}{pressure coefficient}

\
\ 
\newline

\uline{Subscripts}

\ 

\entry{m}{manufactured}
\entry{d}{design intent}
\entry{1}{upstream station}
\entry{2}{downstream station}
\entry{o}{stagnation quantity}
\entry{b}{baseline value}
\end{nomenclature}

%%%%%%%%%%%%%%%%%%%%%%%%%%%%%%%%%%%%%%%%%%%%%%%%%%%%%%%%%%%%%%%%%%%%%%

\section{Introduction}

The geometry of manufactured compressor blades inevitably differs from the nominal, or design intent, geometry specified by the designer due to noise in the manufacturing process or in-service erosion. The impact of geometric variability is characterized by an increase in performance variability and a degradation of the mean performance. For example, Garzon demonstrated that the mean loss coefficient of a flank-milled integrally bladed rotor (IBR) increased by 23\% due to manufacturing variability \cite{garzon_thesis}. Not surprisingly, this shift in performance scales with the level of noise in the manufacturing process, suggesting that tightening manufacturing tolerances can improve the mean performance. Reducing manufacturing tolerances increases the manufacturing costs, so it is important to reduce tolerances only in the regions of the blade that have the largest impact on the mean performance.

An alternative approach to improve the mean performance is to perform robust design optimization of the nominal blade geometry. Robust optimization seeks to find a design whose performance remains relatively unchanged when variability is introduced. Robust optimization is typically probabilistic, seeking to optimize the mean, variance or failure probability of the design \cite{garzon_thesis, kumar_2006_2, bestle_2010}. Most previous work has considered optimizing the nominal blade geometry for a fixed level of manufacturing noise, i.e. fixed manufacturing tolerances. The notable exception is the recent work of Goodhand \textit{et al.}, which demonstrated that the robust blade geometry depends on the level of geometric noise only when a ``switch'' in the dominant loss mechanism occurs on some of the manufactured blades \cite{goodhand_2014}. They provide the example of leading edge flow separation that occurs on some manufactured blades when the leading edge radius of curvature is reduced beyond a critical value, dramatically increasing the blade losses. The optimal leading edge geometry therefore depends on the level of manufacturing noise, illustrating the potential interactions of geometry and tolerance optimization.

In this paper, we present an approach for designing the manufacturing tolerances to improve the mean performance of manufactured blades. The effect of tightening manufacturing tolerances is modelled by reducing the geometric variability. An optimization framework is developed that reduces manufacturing variability in the regions of the blade with the largest impact on the mean performance, and a gradient-based approach for determining the optimal tolerances is presented. This tolerance optimization framework is incorporated into a geometry design framework to perform optimization of the blade geometry and manufacturing tolerances. The interdependence of the optimal blade geometry and manufacturing tolerances is analyzed for an exit guide vane, and recommendations for the design process are made.

\section{Random Field Model of Manufacturing Variability}

Previous studies of manufactured compressor blades have shown that the surface deviations are often normally distributed \cite{beachkofski_2004, duffner_thesis, garzon_thesis, sinha_2008}. We therefore model the manufacturing deviation as a Gaussian random field $e(s,\samp)$ where $s$ indexes the spatial location on the blade and $\samp$ indexes the probability space. The random field $e(s,\samp)$ describes the error between the manufactured surface and the nominal surface in the normal direction at the point $s$ on the nominal blade surface:

\begin{equation}
x_m(s,\samp) = x_d(s) + e(s,\samp)\hat{n}(s).
\end{equation}
This random field $e(s,\samp)$ is uniquely defined by its mean $\overline{e}(s) := \E[e(s,\samp)]$ and covariance function $C(s,s') := \E[(e(s,\samp) - \overline{e}(s))(e(s',\samp) - \overline{e}(s'))]$. 

A scaling study performed by Garzon showed that for realistic levels of manufacturing variability, the mean shift in the loss coefficient resulting from changes in the mean geometry account for less than 15\% of the total mean shift caused by manufacturing variations \cite{garzon_thesis}. Since the performance impact resulting from changes in the mean geometry are small compared to the impact associated with increased variance, changes in the mean geometry are not modeled and the error field is assumed to have zero mean. Previous studies of manufacturing variations have also shown that the correlation length (the characteristic length scale over which the correlation function $\rho(s,s') := C(s,s')/\sigma(s)\sigma(s')$ decays) is smaller near the leading edge of compressor blades \cite{garzon_thesis, lange_2012}. This results in variability in the bluntness of the leading and trailing edges. We simulate the manufacturing errors using the squared exponential correlation function, which produces manufacturing errors with continuous curvature:

\begin{equation}
\rho(s,s') = \exp\left(-\frac{|s - s'|^2}{2L^2}\right).
\end{equation}
The correlation length $L := [\tilde{L}(s)\tilde{L}(s')]^{1/2}$ is chosen to be smaller near the leading edge to produce variations in the leading edge bluntness:

\begin{equation}
\tilde{L}(s) = L_0 + (L_{LE} - L_0)\exp(-s^2/w^2),
\end{equation}
where $L_{LE} < L_0$, and the width parameter $w$ is equal to the leading edge radius.

\subsection{The Karhunen-Lo\`{e}ve Expansion}

The manufacturing error field $e(s,\samp)$ is simulated using Karhunen-Lo\`{e}ve (K-L) expansion. The K-L expansion, also referred to as the proper orthogonal decomposition (POD), represents a random field as a spectral decomposition of its covariance function:
\begin{equation}
e(s,\samp) = \bar{e}(s) + \sum_{i \geq 1} \sqrt{\lambda_i} \phi_i(s) \xi_i(\samp),
\label{eq:kl_full}
\end{equation}
where $\lambda_i$ and $\phi_i(s)$ are the eigenvalues and eigenfunctions of the covariance function, respectively, and the $\xi_i(\samp)$ are mutually uncorrelated random variables with zero mean and unit variance \cite{loeve}. For a Gaussian random field, the $\xi_i(\samp)$ are i.i.d. standard normal random variables \cite{lemaitre}.

\section{Tolerance and Geometry Design}

Tolerance design changes the allowable level of geometric variability in the population of manufactured blades. In general, the optimal tolerance design depends on the nominal blade geometry since changing the nominal geometry changes the location of flow features whose performance impact is sensitive to geometric variability, i.e. the location of transition or flow separation. The optimal blade geometry is found to depend on the tolerances if changing the tolerances results in a switch in the dominant loss mechanism for some of the manufactured blades. We now formulate the tolerance and geometry design optimization problems, and illustrate the interdependence of tolerance and geometry design optimization for an exit guide vane geometry.

\subsection{Tolerance Design}

The precision of a manufacturing process can be specified in terms of its spread, defined as some multiple of the standard deviation of the dimension of interest, where a dimension refers to the size of a particular feature \cite{kane_1986}. To specify manufacturing tolerances around a compressor blade surface, we specify the process spread, equal to the standard deviation $\sigma(s)$, at every point on the blade surface. 
%The standard deviation field $\sigma(s)$ is represented using a cubic B-spline basis:
%\begin{equation}
%\sigma(s) = \sum_{i = 1}^{N_\sigma}\sigma_i B_{i}(s),
%\end{equation}
%where the $B_{i}(s)$ are cubic B-spline basis functions. This choice of basis ensures that the manufactured blade surfaces have continuous curvature, i.e. have no sharp edges. The manufacturing tolerances are parameterized by the vector $\tol = \{\sigma_i\}_{i=1}^{N_\sigma}$.

In the proposed tolerance design framework, the process spread $\sigma(s)$ is optimized to minimize the mean loss. Since the loss coefficient depends on the flow incidence, we design the manufacturing tolerances to minimize the mean loss coefficient over a range of incidence angles by adopting a multi-point design strategy. A set of $N_p$ design points is chosen, each corresponding to an incidence angle. An objective function is defined as the weighted sum of the loss coefficient at each incidence angle, i.e.
\begin{equation}
J := \sum_{i=1}^{N_p} w_i \loss(\alpha_i),
\end{equation}
where the weight vector $\mathbf{w} = \{w_1,...,w_{N_p}\}$ weights the relative importance of different incidence angles. 

%The weight vector is chosen so that the multi-point objective function $J$ approximates a weighted integral of the loss coefficient:
%\begin{equation}
%\sum_{i=1}^{N_p} w_i \loss(\beta_i) \approx \int_\mathcal{B} W(\beta) \loss(\beta)\, d\beta,
%\end{equation}
%where $W(\beta)$ is a weighting function and $\mathcal{B} = [\beta_1^l,\beta_1^u]$ is the incidence range of interest. Choosing $N_p = 3$ equispaced points and a loss weight vector $\mathbf{w} = \{1/4, 1/2, 1/4\}$ approximates the average loss computed using a trapezoidal integration rule.

The tolerance design process determines where on the blade surface the tolerances should be reduced to have the most benefit, i.e. to provide the greatest reduction in the mean loss coefficient. This is achieved by constraining the level of manufacturing variability, and minimizing the weighted mean loss $\E[J]$. To quantify the level of variability, we introduce the function $V(\tol)$, equal to the integral of the standard deviation of the random field that models the manufacturing variations over the blade surface:
\begin{equation}
V(\tol) = \int_S \sigma(s)\, ds.
\label{eq:variability}
\end{equation}
Reducing $V$ implies stricter tolerances and therefore increased manufacturing cost. In the optimization, the standard deviation is also bounded from above by a specified value $\sigma_{max}$ to ensure that the optimizer does not trade increases in the variability in regions of low mean loss sensitivity for decreases in high sensitivity regions and drive $\sigma(s)$ to zero. The optimization problem for designing the tolerances is summarized below:
\begin{equation}
\begin{aligned}
& \tol^* =
& & \underset{\tol}{\arg \min}
& & \E[J(\samp;\tol)] \\
& & & \ \ \ \ \text{s.t.}
& & V(\tol) = V_b \\
& & & & & 0 \leq \sigma(s) \leq \sigma_{max}.
\end{aligned}
\label{eq:tol_opt}
\end{equation}

\subsection{Geometry Design}

The blade performance can also be improved by optimizing the nominal blade geometry. Geometry optimization can be either deterministic, where the performance in the absence of geometric variability is optimized, or robust, where a performance statistic in the presence of geometric variability is optimized. The deterministic approach can be summarized by the following optimization statement, where the weighted loss is minimized subject to a constraint on the flow turning $\turn$:
\begin{equation}
\begin{aligned}
& \des_{det}^* =
& & \underset{\des}{\arg \min}
& & J(\des) \\
& & & \ \ \ \ \text{s.t.}
& & \turn(\des) = \turn_b.
\end{aligned}
\label{eq:det_opt}
\end{equation}
The robust equivalent of (\ref{eq:det_opt}) replaces the weighted loss and turning by their mean values:
\begin{equation}
\begin{aligned}
& \des_{rob}^* =
& & \underset{\des}{\arg \min}
& & \E[J(\samp; \des)] \\
& & & \ \ \ \ \text{s.t.}
& & \E[\turn(\samp; \des)] = \turn_b.
\end{aligned}
\label{eq:rob_opt}
\end{equation}

Solving the robust optimization problem (\ref{eq:rob_opt}) incurs a higher computational cost than solving the deterministic optimization problem (\ref{eq:det_opt}) since estimating the performance statistics requires multiple evaluations of the performance. If the optimal geometry is insensitive to geometric variability, then $\des^*_{det} = \des^*_{rob}$ and the computationally cheaper deterministic design approach can be used. In this case, the geometry and tolerance optimizations are decoupled, and the tolerances can be optimized after the nominal blade geometry by solving (\ref{eq:tol_opt}).

\subsection{Simultaneous Tolerance and Geometry Design}

When the optimal geometry does depend on the level of geometric variability, optimizing the blade geometry first will result in a sub-optimal design when the level of variability is changed by the tolerance optimization step. In this case, the geometry and tolerance optimizations are coupled and should be performed simultaneously. The simultaneous geometry and tolerance optimization is formulated by combining the tolerance and robust geometry optimization problems given by (\ref{eq:tol_opt}) and (\ref{eq:rob_opt}):
\begin{equation}
\begin{aligned}
& (\des^*,\tol^*) =
& & \underset{\des,\tol}{\arg \min}
& & \mathbb{E}[J(\samp; \des,\tol)] \\
& & & \ \ \ \ \text{s.t.}
& & \E[\turn(\samp; \des, \tol)] = \turn_b \\
& & & & & V(\tol) = V_b \\
& & & & & 0 \leq \sigma(s) \leq \sigma_{max}.
\end{aligned}
\label{eq:sim_opt}
\end{equation}

The above optimization problems are solved using a gradient-based approach. For the optimizations involving the mean performance, the sample average approximation (SAA) method is used to compute an approximate solutions. The SAA method replaces the mean loss and mean turning with Monte Carlo estimates computed using a fixed set of random inputs, i.e. a fixed set of realizations $\{\bxi^{(n)}\}_{n=1}^{N}$ of the random input vector used to generate realizations of the error field $e(s,\samp)$ \cite{rubinstein_shapiro}. Fixing the random inputs transforms the stochastic optimization problems into deterministic optimization problems. The deterministic optimization problems are solved using the sequential quadratic programming (SQP) method, which is a gradient-based method for solving nonlinear constrained optimization problems \cite{powell_1978}.

\subsection{Application}

We optimize the nominal geometry and tolerances of a two-dimensional fan exit guide vane developed by UTRC with design inlet Mach number $M_1 = 0.73$, Reynolds number $1.0 \times 10^6$, and axial velocity density ratio of 1.12 \cite{stephens_hobbs_1979}. A turbulence intensity of 4\% was selected to reflect conditions in a typical compressor stage \cite{camp_1995}. All flow solutions are computed using the Multiple Blade Interacting Streamtube Euler Solver (MISES) code, which includes models for both natural and bypass transition, and can also model shock-induced and trailing edge separation \cite{mises_manual}. 

The nominal blade geometry is parameterized using a combination of Chebyshev polynomial modes, acting in the direction normal to the baseline blade surface, and a stagger angle mode to rotate the blade around the leading edge. Five Chebyshev modes were used, giving a total of $N_d = 6$ nominal geometry design parameters. The blade thickness is constrained by applying equal and opposite perturbations to corresponding points on the pressure and suction surfaces so that only the camber and stagger angle are modified. 

The standard deviation field $\sigma(s)$ is represented using a cubic B-spline basis. This choice of basis ensures that the manufactured blade surfaces have continuous curvature, i.e. have no sharp edges. The standard deviation is represented using $N_\sigma = 41$  basis functions, and the knot spacing is reduced near the leading edge, since previous studies of the impact of geometric variability on compressor performance suggest that most of the increase in loss results from imperfections near the leading edge \cite{garzon_2003}.

\subsubsection{Single-point Optimization}

We first perform a single-point optimization to minimize the loss at the design incidence. Deterministic (Det) and robust (Rob) optimizations of the geometry are performed, as well as a simultaneous (Sim) optimization of the geometry and tolerances. The tolerances of the deterministic and robust optimal blade geometries are also optimized. In all tolerance optimizations, the standard deviation is constrained from above by $\sigma_b/c = 8.0 \times 10^{-4}$, roughly twice the level of variability observed in point-milled compressor blades \cite{garzon_thesis}. The total variability $V$ is constrained to be 98\% of the value computed by evaluating (\ref{eq:variability}) with $\sigma/c = 8.0 \times 10^{-4}$.

Table \ref{tab:uti_sp} summarizes the performance of the single-point optimized designs. The second row gives the mean loss calculated using uniform manufacturing variability ($\sigma_b/c = 8.0 \times 10^{-4}$), and the last row gives the mean loss calculated using optimized tolerances. In the absence of variability, the deterministic optimum has the lowest loss coefficient, 11\% lower than the loss of the baseline design. The robust geometry optimization reduces the mean loss by 8\%, and has the lowest mean loss at the baseline level of manufacturing variability. The simultaneous optimization has lower mean loss than the other designs with optimized tolerances, and reduces the mean loss by 11\% relative to the baseline geometry with uniform tolerances.

\begin{table}[htbp]
\begin{center}
\begin{tabular}{|c|c|c|c|}
\hline
 & $\loss$ & $\E[\loss(\tol_b)]$ & $\E[\loss(\tol^*)]$ \\ \hline
Base & $2.24 \times 10^{-2}$ & $2.27 \times 10^{-2}$ & $2.25 \times 10^{-2}$ \\
\hline
Det & $2.00 \times 10^{-2}$ & $2.14 \times 10^{-2}$ & $2.06 \times 10^{-2}$ \\
\hline
Rob & $2.06 \times 10^{-2}$ & $2.09 \times 10^{-2}$ & $2.06 \times 10^{-2}$ \\
\hline
Sim & $2.03 \times 10^{-2}$ & $2.10 \times 10^{-2}$ & $2.05 \times 10^{-2}$ \\
\hline
\end{tabular}
\end{center}
\caption{Performance comparison of the baseline UTRC, single-point deterministic optimal, single-point robust optimal, and single-point simultaneous optimal designs.}
\label{tab:uti_sp}
\end{table}

\begin{figure*}[htbp]
\centering
\subfigure
{\includegraphics[width=0.49\textwidth]{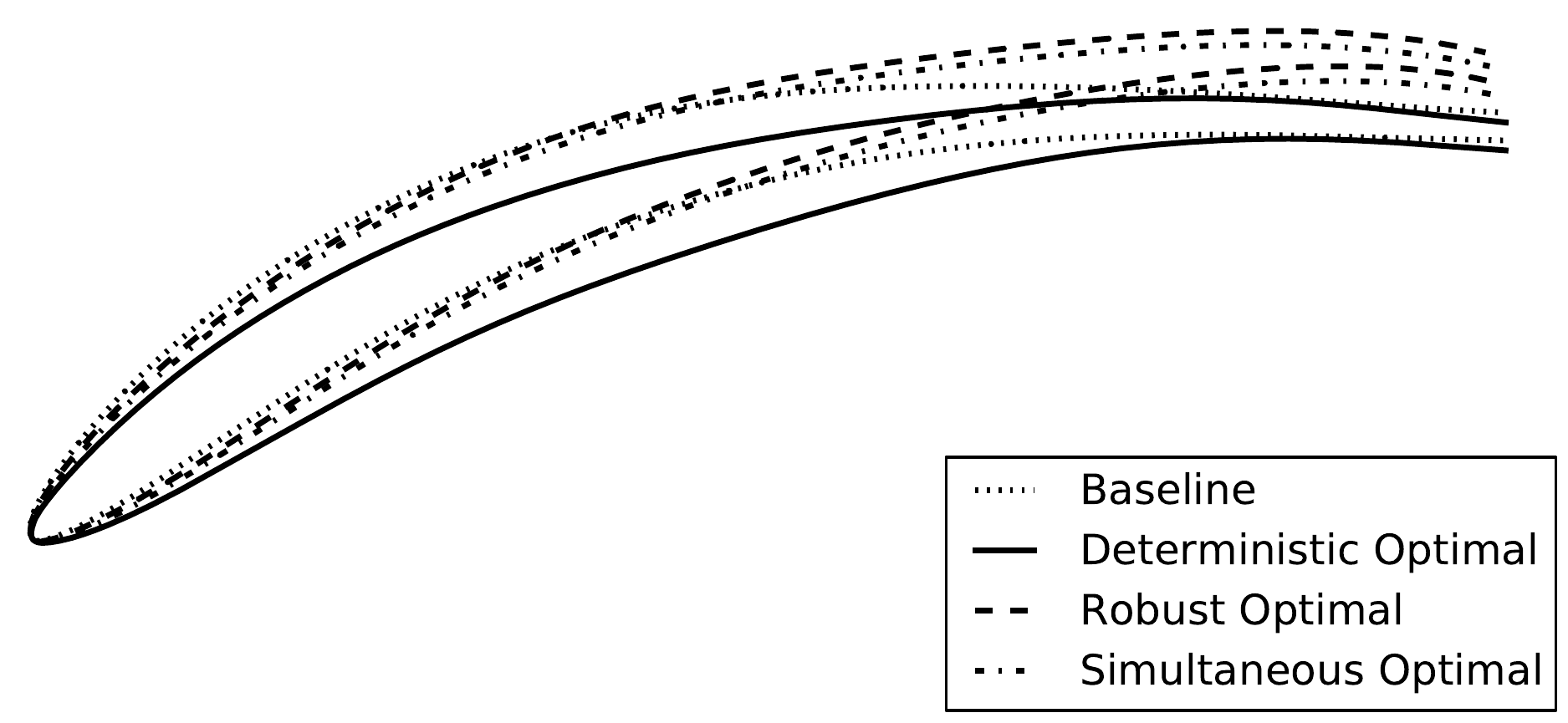}}
\subfigure
{\includegraphics[width=0.49\textwidth]{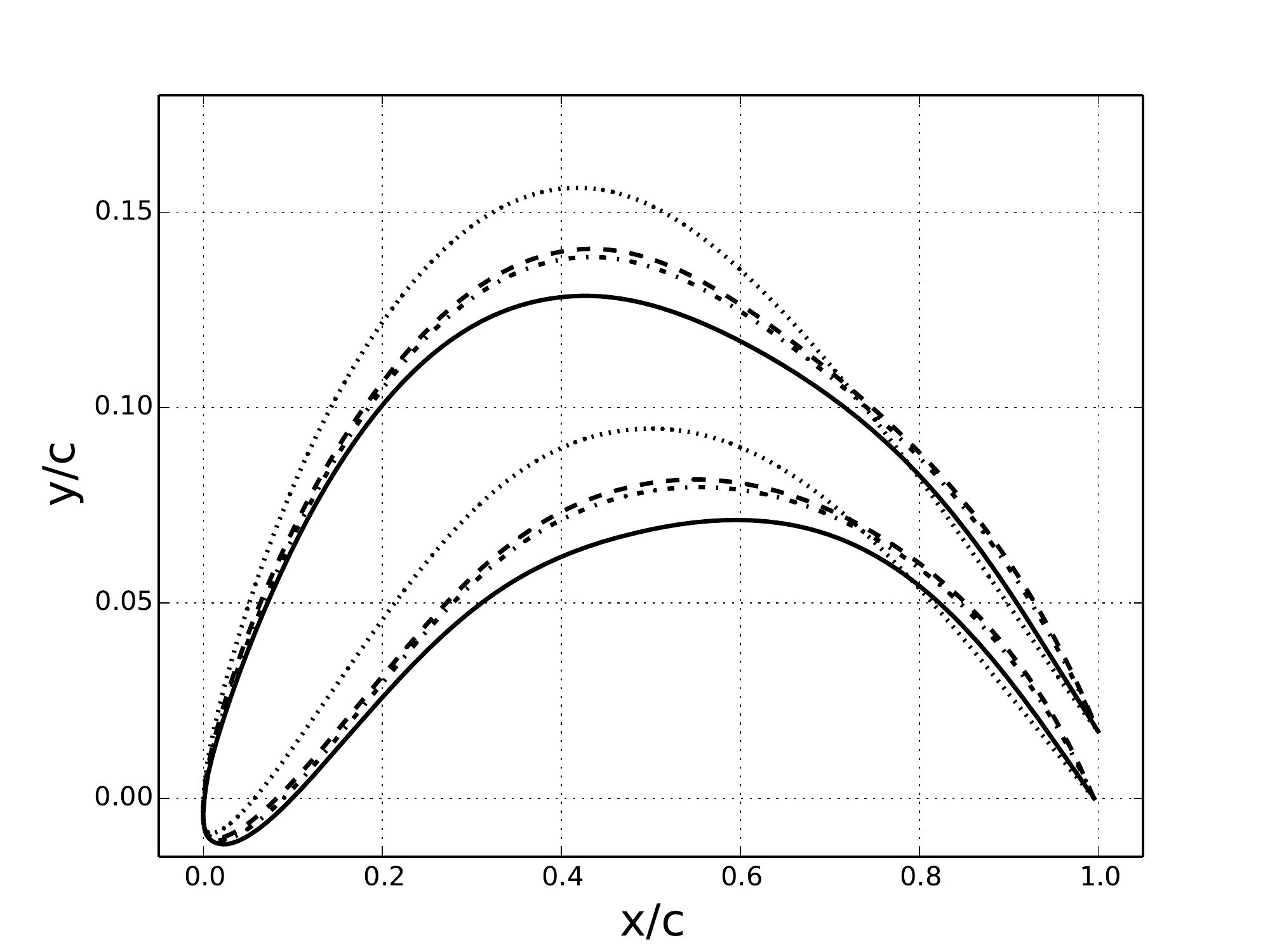}}
\caption{Single-point optimal redesigned UTRC blades.}
\label{fig:uti_sp_blades}
\end{figure*}

\begin{figure}[htbp]
\centering
\includegraphics[width=0.5\textwidth]{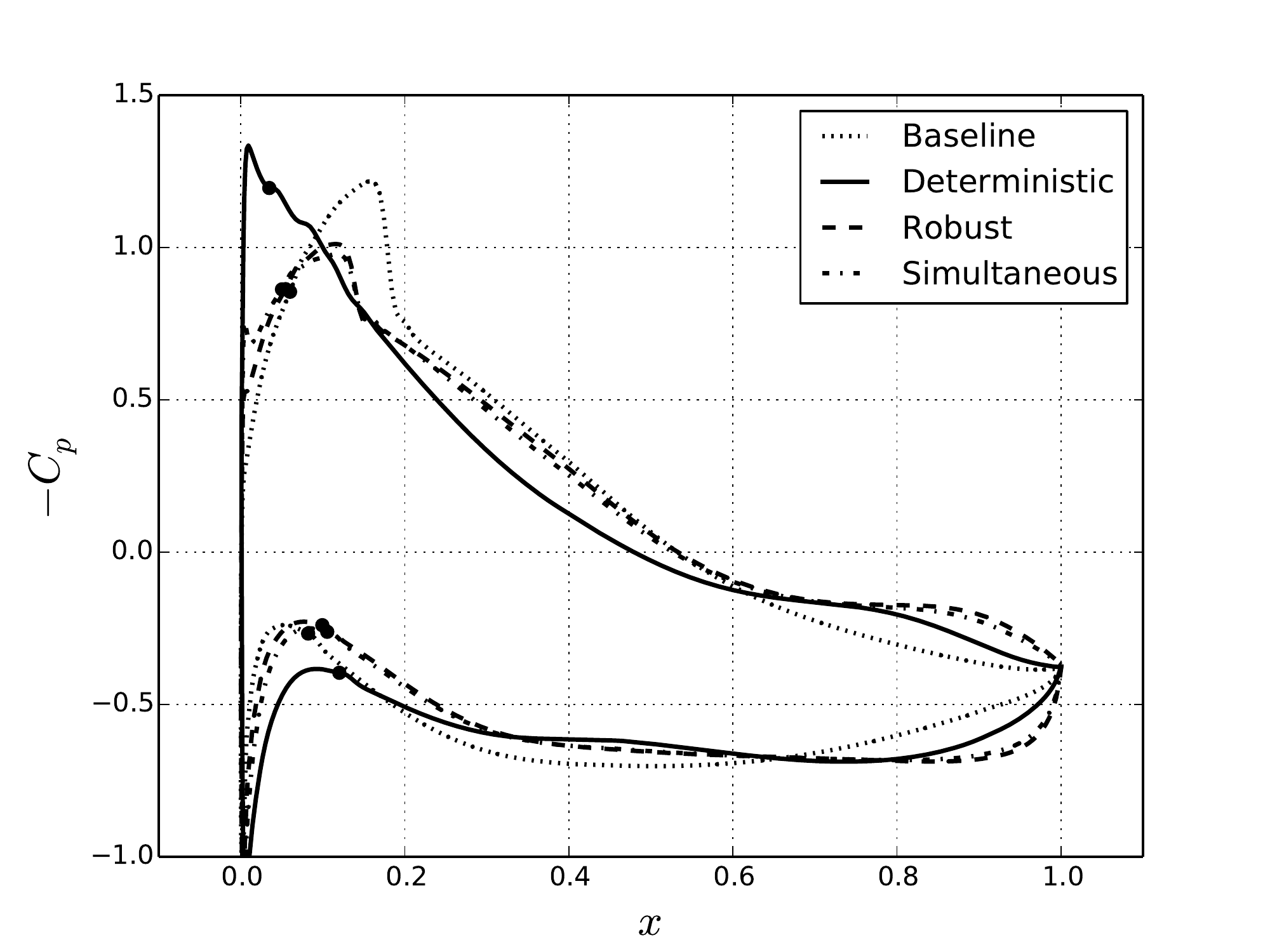}
\caption{Pressure coefficient profiles for the baseline and single-point optimized UTRC blades. The dots denote the location of transition.}
\label{fig:uti_cp_comb_sp}
\end{figure}

The optimal blade geometries are shown in Figure \ref{fig:uti_sp_blades}, and Figure \ref{fig:uti_cp_comb_sp} compares the pressure coefficient profiles for the three nominal blade designs. The deterministic optimization eliminates the shock by reducing the camber of the first 75\% of the chord, which reduces the losses generated on the suction side of the blade. The pressure profiles for the robust and simultaneous optimal blades are similar, with a shock forming on the suction side for both nominal geometries. The peak Mach number on the suction side of the robust and simultaneous optimal blades are 1.15 and 1.14, respectively, compared to 1.32 for the deterministic optimal blade and 1.26 for the baseline blade. The reduction in shock strength reduces the shape factor downstream of the shock relative to the baseline geometry, resulting in lower boundary layer losses.

Deterministic optimization produces a blade that is not robust to manufacturing variability: the mean loss of the manufactured blades is 7\% higher than the nominal loss, and is higher than the mean loss of the other two optimized blade geometries. The deterministic optimization eliminates the suction side shock, but also increases the peak suction side Mach number relative to the baseline design. Some of the manufactured blade realizations lie in the region of the design space where a switch in dominant loss mechanism occurs. Shocks form on the suction side of the blade near the leading edge for some of the manufactured blades, and shock induced separation occurs for 10\% of the manufactured blades. 

The robust design avoids the switch in flow mechanisms by reducing the peak suction side Mach number. Introducing manufacturing variations results in shock induced flow separation on 5\% of the blades, resulting in a 2\% increase in mean loss. Simultaneous optimization also avoids the switch in flow mechanisms by changing the geometry and reducing the manufacturing variability, and the flow remains attached for all manufactured blades.

The single-point optimal geometry clearly depends on the level of geometric variability, and thus the manufacturing tolerances, since deterministic and robust optimization produce different designs. The dependence of the optimal geometry on the level of variability is a result of a switch in the dominant loss mechanism that occurs on some of the manufactured blades. 

\begin{figure*}[htbp]
\centering
\subfigure[Deterministic Optimal Tolerances]{\includegraphics[width=0.49\textwidth]{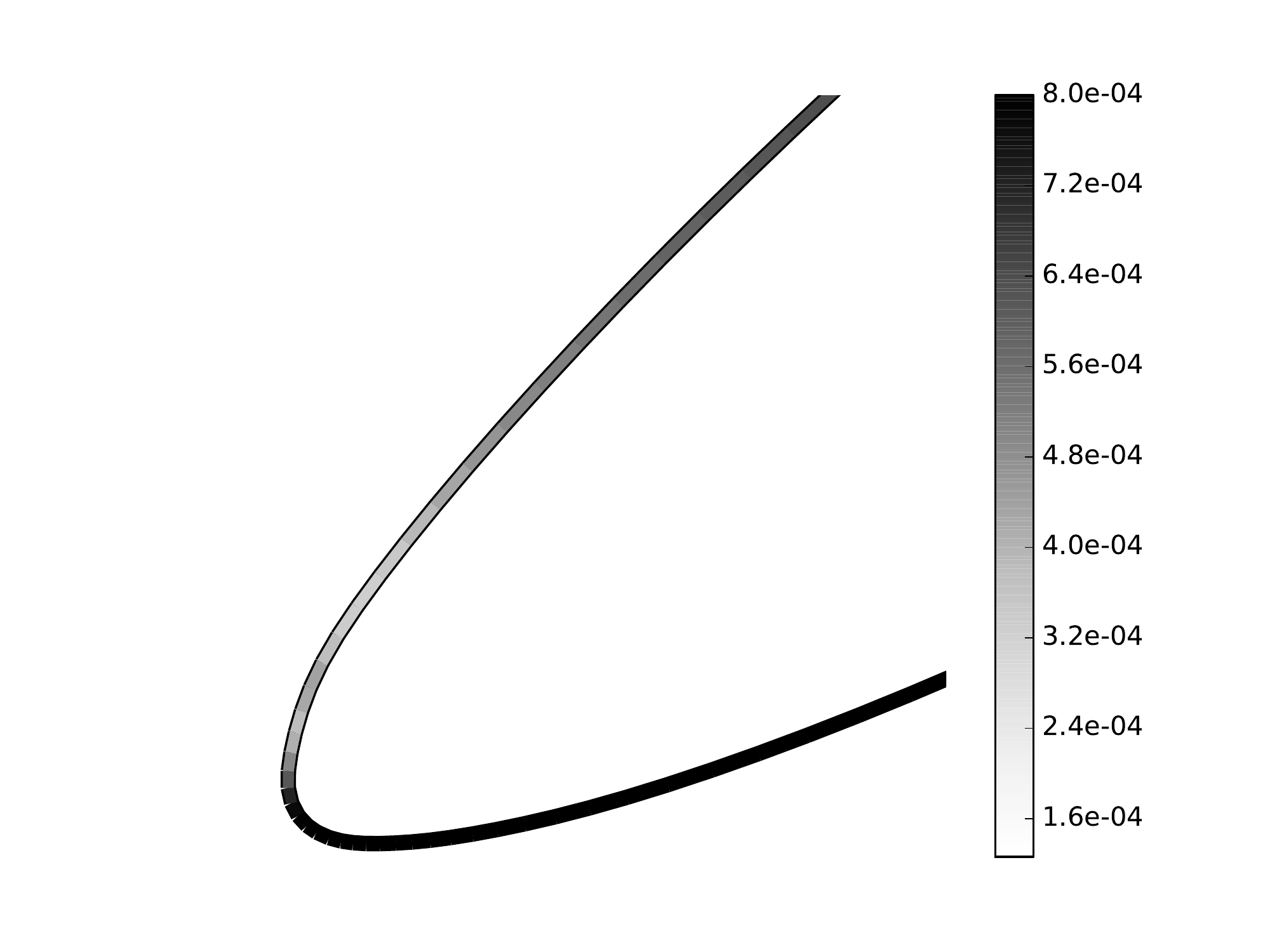}}
\subfigure[Robust Optimal Tolerances]{\includegraphics[width=0.49\textwidth]{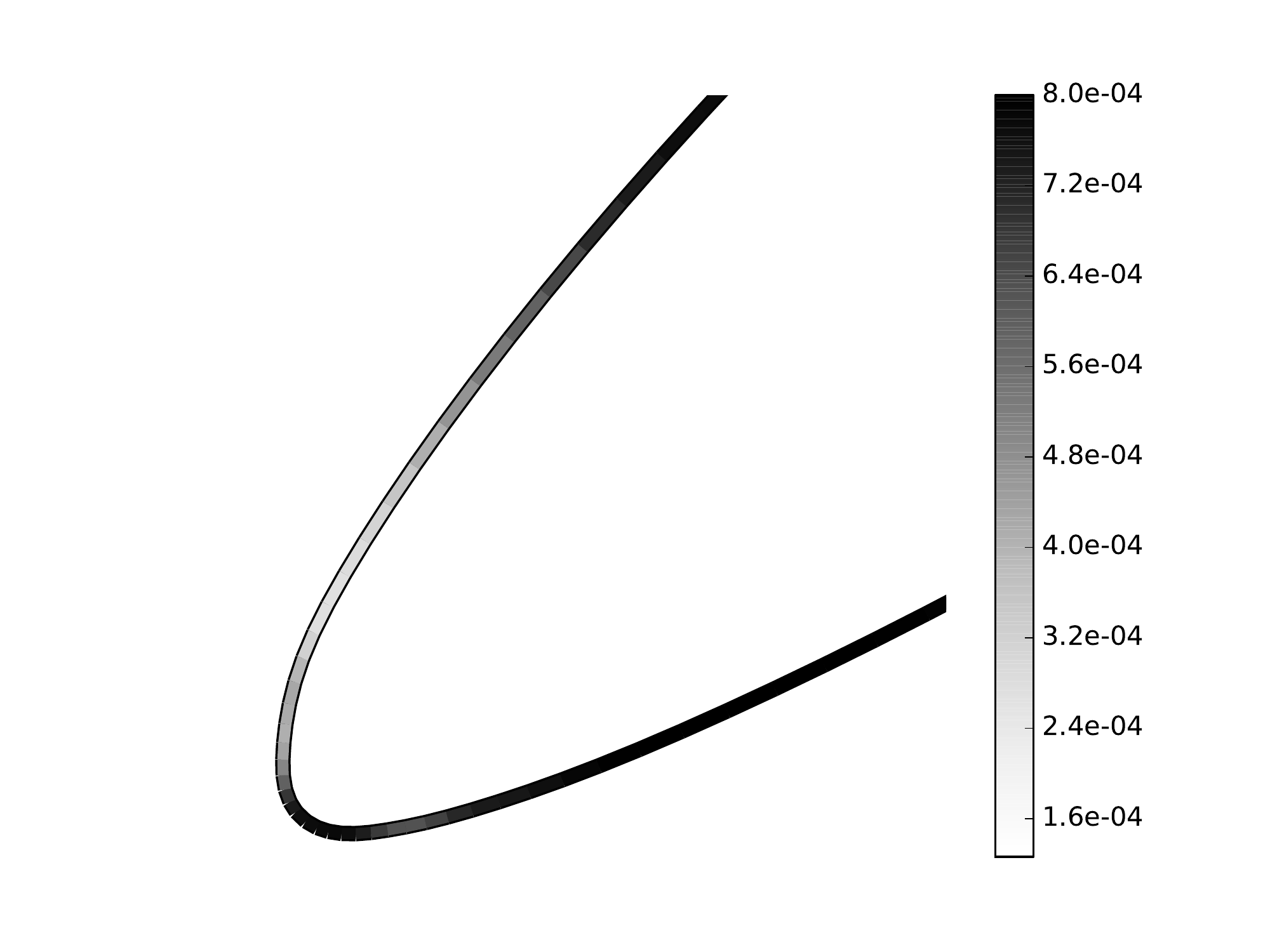}}
\subfigure[Simultaneous Optimal Tolerances]{\includegraphics[width=0.49\textwidth]{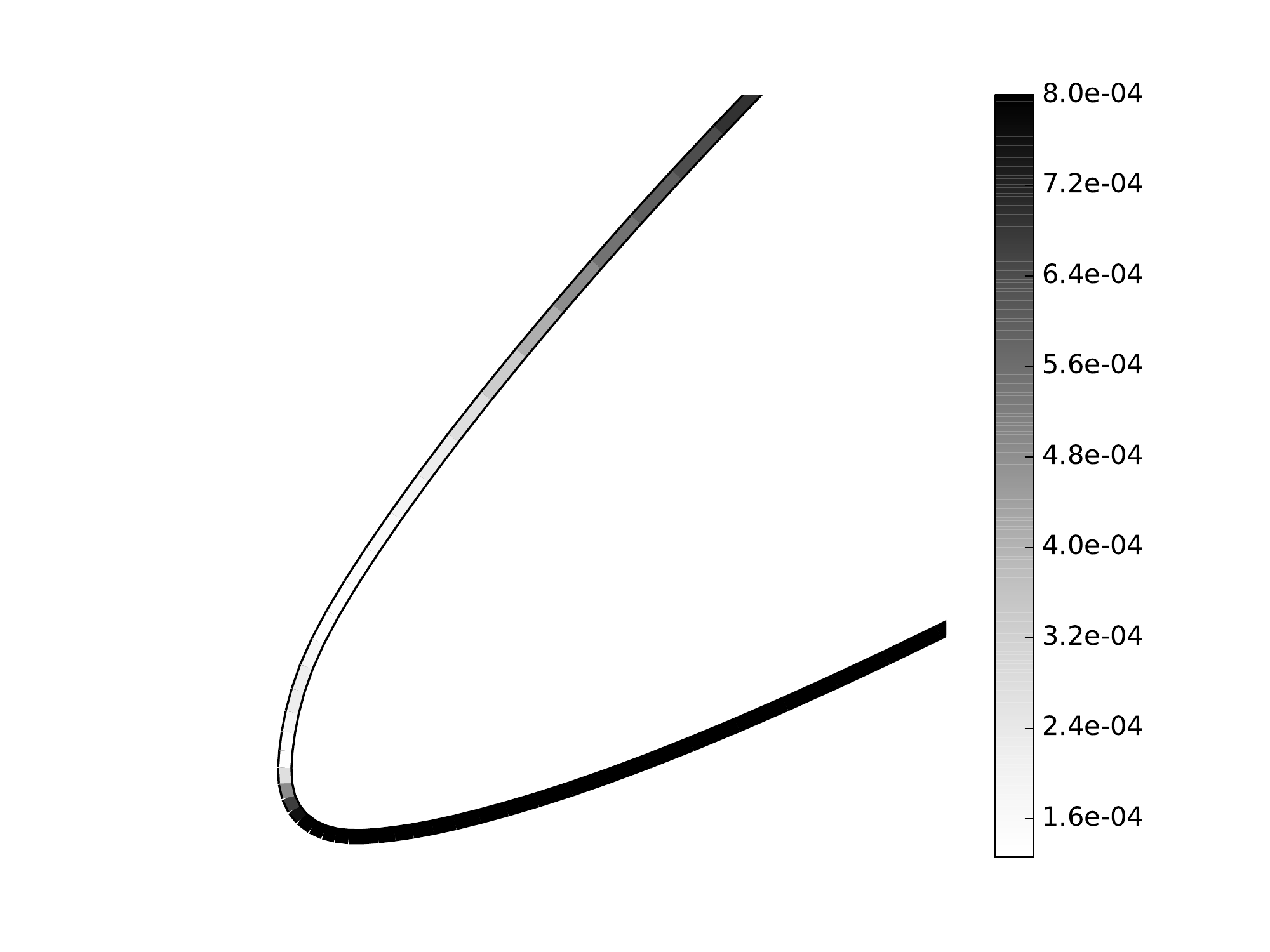}}
\caption{Optimal standard deviation $\sigma(s)/c$ for the single-point optimized UTRC blades. The lower surface is the pressure side, and the upper surface is the suction side.}
\label{fig:uti_sp_opt_tol}
\end{figure*}

The optimal tolerances are compared in Figure \ref{fig:uti_sp_opt_tol}. We only show the leading edge since the optimal standard deviation outside this region is equal to the baseline value of $\sigma_b/c = 8.0 \times 10^{-4}$. The reduction in variability on the deterministic optimal blade extends further down the suction side of the blade than for the other geometries. This reduces the percentage of blades with suction side flow separation to 3\%. The optimal tolerances for the robust optimal blade are distributed according to a ``double bow-tie'' pattern: the tolerances are reduced on either side of the leading edge. The reduction in variability on the pressure side is small relative to that on the suction side. The pressure side tolerances ensure flow separation does not occur on any manufactured robust optimal blades. To quantify the error between two tolerance schemes $\sigma_i(s)$ and $\sigma_j(s)$, we compute the integrated error:
\begin{equation}
e_{i,j} = \frac{\displaystyle\int_S \bigl\lvert[\sigma_{max} - \sigma_i(s)] - [\sigma_{max} - \sigma_j(s)]\bigr\rvert \ ds}{\displaystyle\int_S [\sigma_{max} - \sigma_i(s)] + [\sigma_{max} - \sigma_j(s)]\ ds}.
\end{equation}
The largest error between the three optimal tolerance schemes is 37\%, illustrating the sensitivity of the optimal tolerances on the nominal blade geometry.

\subsubsection{Multi-point Optimization}

\begin{figure*}[htbp]
\centering
\subfigure
{\includegraphics[width=0.49\textwidth]{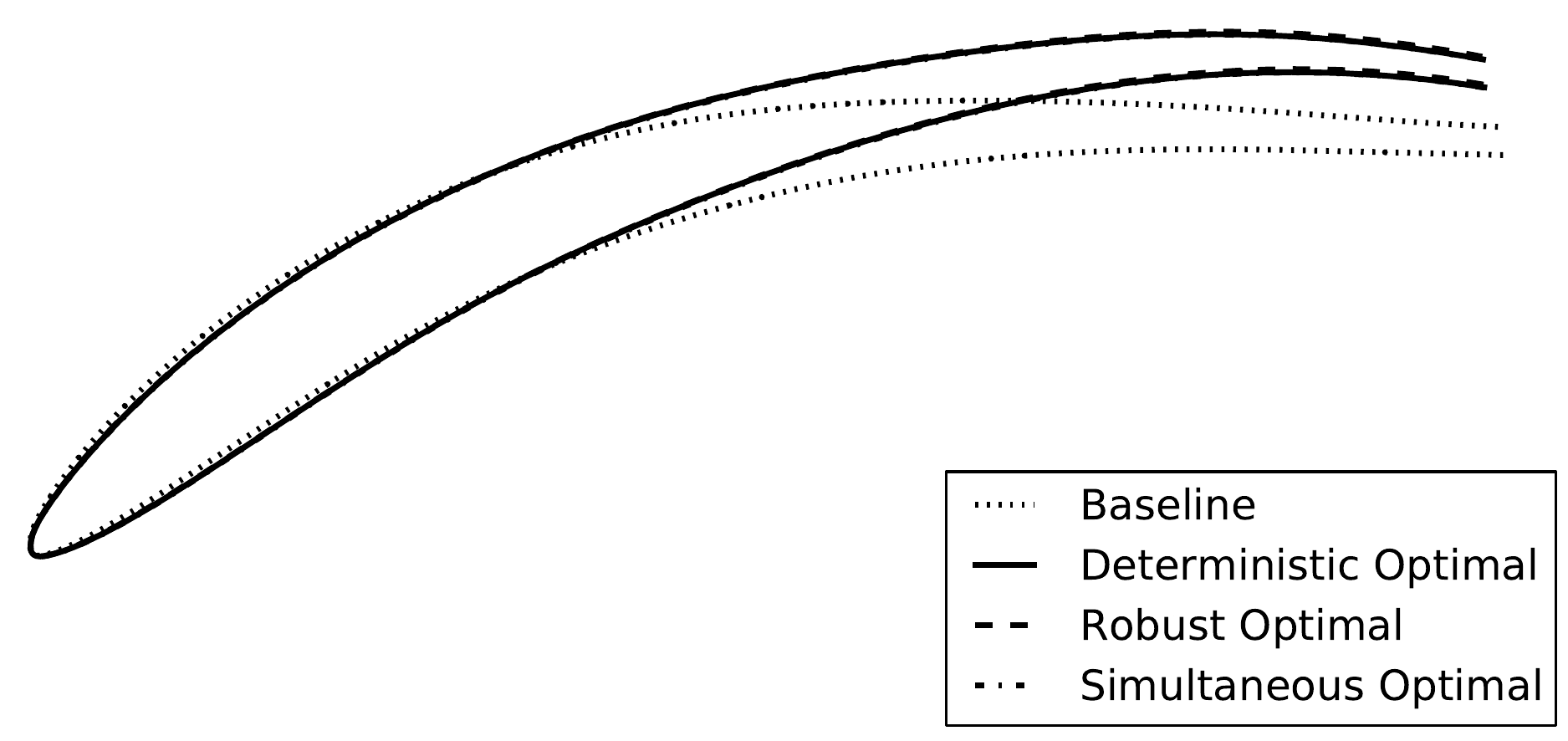}}
\subfigure
{\includegraphics[width=0.49\textwidth]{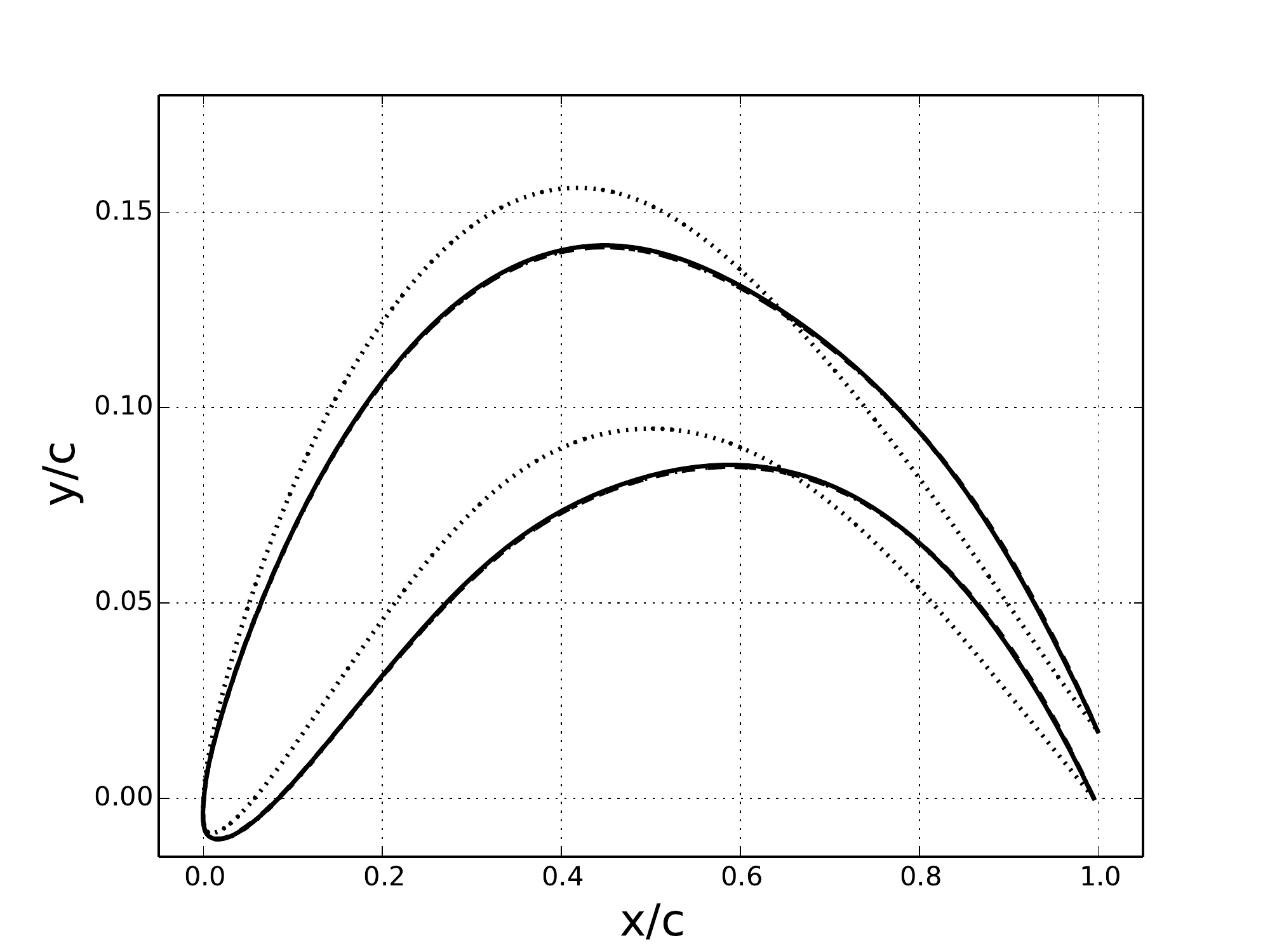}}
\caption{Multi-point optimal redesigned UTRC blades.}
\label{fig:uti_mp_blades}
\end{figure*}

\begin{figure*}[htbp]
\centering
\subfigure[Uniform Tolerances]{\includegraphics[width=0.49\textwidth]{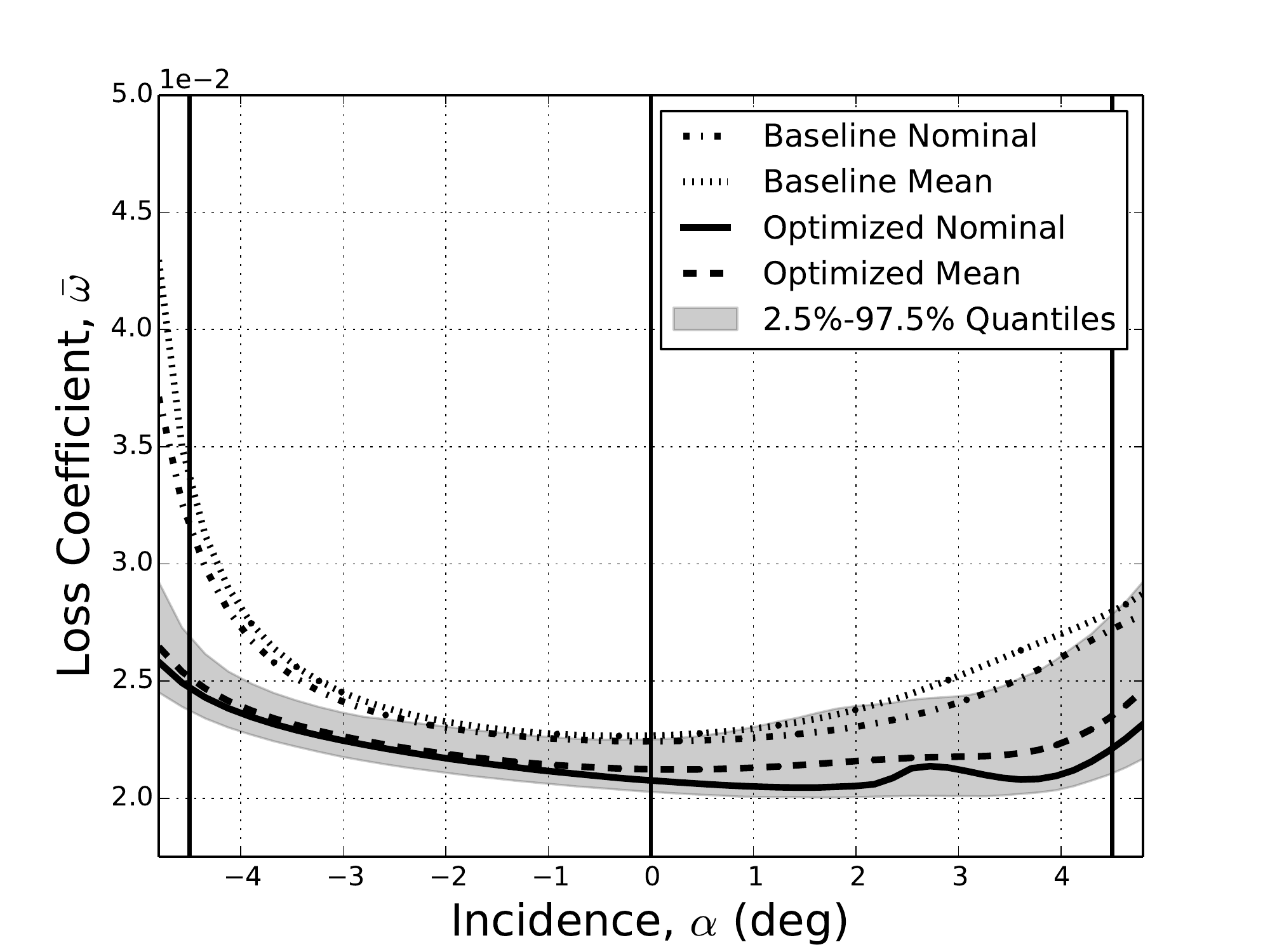}}
\subfigure[Optimized Tolerances]{\includegraphics[width=0.49\textwidth]{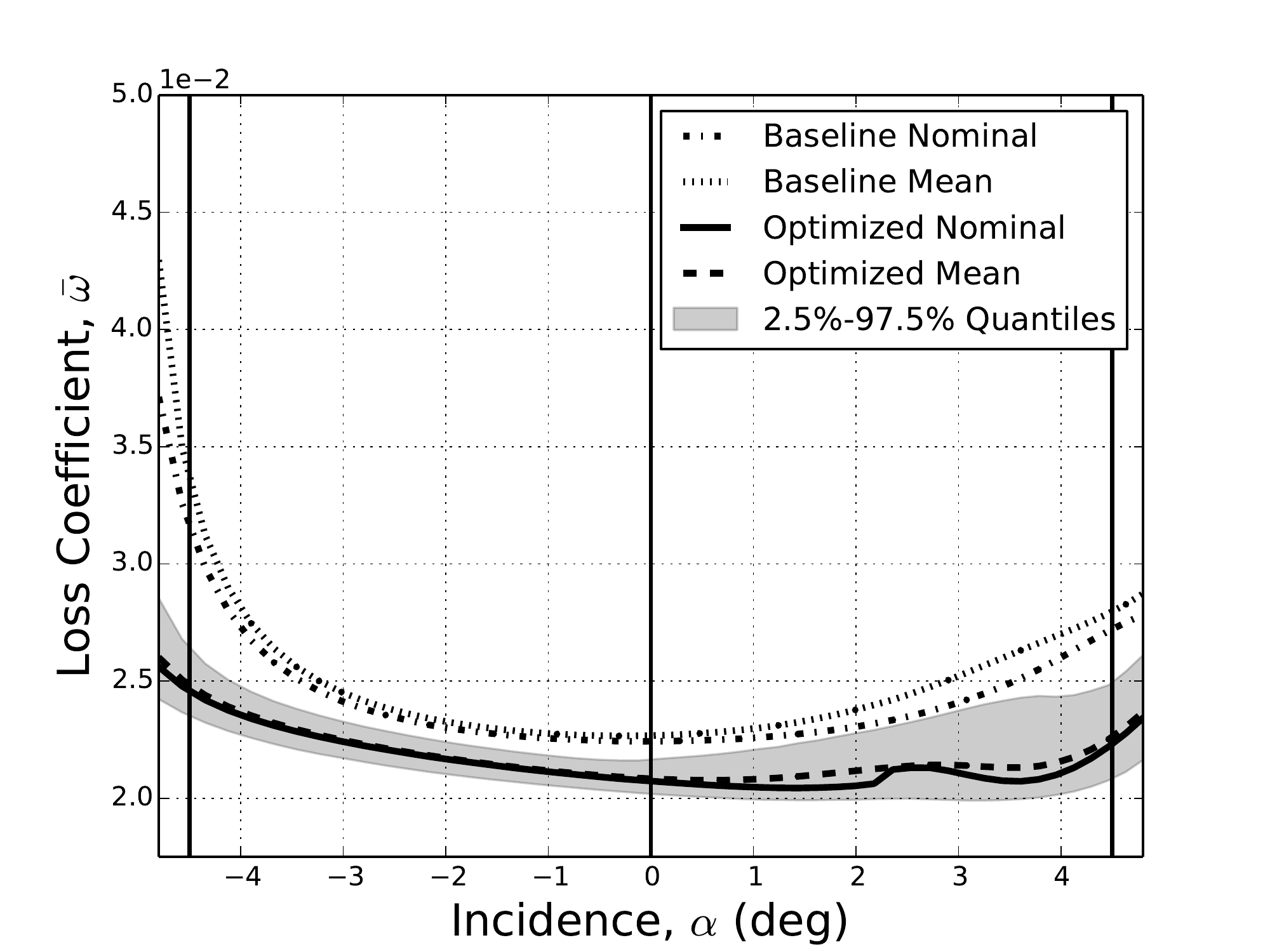}}
\caption{Loss buckets for the multi-point optimized UTRC blade.}
\label{fig:uti_mp_loss_buckets}
\end{figure*}

Multi-point optimizations were performed for the UTRC blade by choosing $N_p = 3$ design points at $\alpha = -4.5^\circ$, $0^\circ$, and $4.5^\circ$ and a loss weight vector $\mathbf{w} = \{1/4,\ 1/2,\ 1/4\}$ so that the weighted loss objective function $J$ approximates the average loss over the incidence range $[-4.5^\circ, 4.5^\circ]$. The multi-point optimal geometries are shown in Figure \ref{fig:uti_mp_blades}. The surfaces of the optimized geometries are nearly identical, and the maximum error between them is 0.5\% of the chord length. The multi-point optimal geometries are similar to the single-point robust geometry: the optimizer reduces the camber of the blade, and decreases the blade stagger. The reductions in the objective function achieved by the various optimizations are summarized in Table \ref{tab:uti_mp}, and a comparison of the loss buckets for the baseline and optimized blades is shown in Figure \ref{fig:uti_mp_loss_buckets}. 

\begin{table}[htbp]
\begin{center}
\begin{tabular}{|c|c|c|c|}
\hline
 & $\sum w_i\loss_i$ & $\sum w_i\mathbb{E}[\loss_i(\tol_b)]$ & $\sum w_i\mathbb{E}[\loss_i(\tol^*)]$ \\ \hline
Base & $2.60 \times 10^{-2}$ & $2.68 \times 10^{-2}$ & $2.64 \times 10^{-2}$ \\
\hline
Det & $2.22 \times 10^{-2}$ & $2.29 \times 10^{-2}$ & $2.23 \times 10^{-2}$ \\
\hline
Rob & $2.22 \times 10^{-2}$ & $2.29 \times 10^{-2}$ & $2.23 \times 10^{-2}$ \\
\hline
Sim & $2.22 \times 10^{-2}$ & $2.29 \times 10^{-2}$ & $2.23 \times 10^{-2}$ \\
\hline
\end{tabular}
\end{center}
\caption{Performance comparison of the baseline UTRC, multi-point deterministic optimal, multi-point robust optimal, and multi-point simultaneous optimal geometries.}
\label{tab:uti_mp}
\end{table}

\begin{figure*}[htbp]
\centering
\subfigure[$\alpha = -4.5^\circ$]{\includegraphics[width=0.49\textwidth]{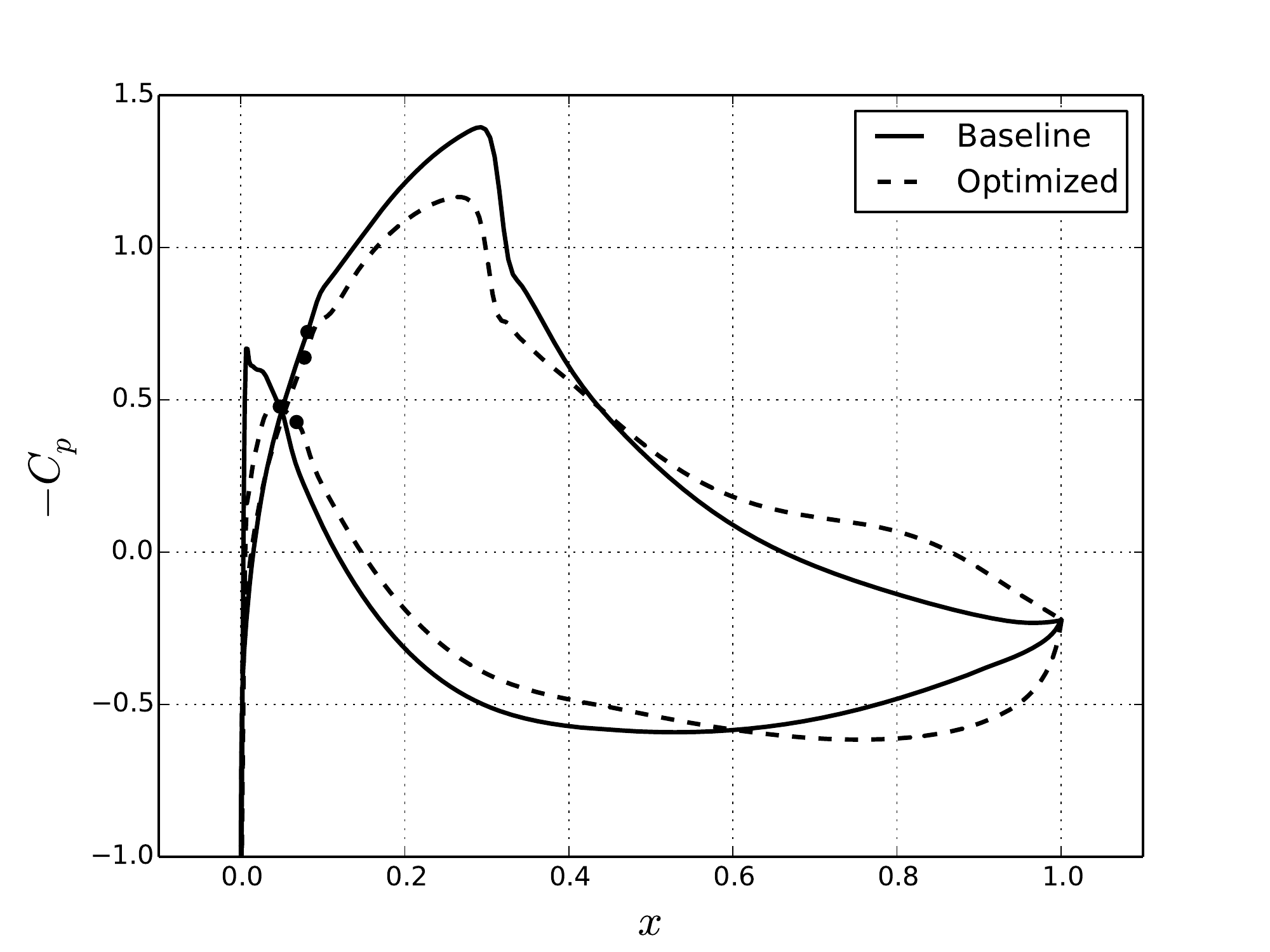}}
\subfigure[$\alpha = 0^\circ$]{\includegraphics[width=0.49\textwidth]{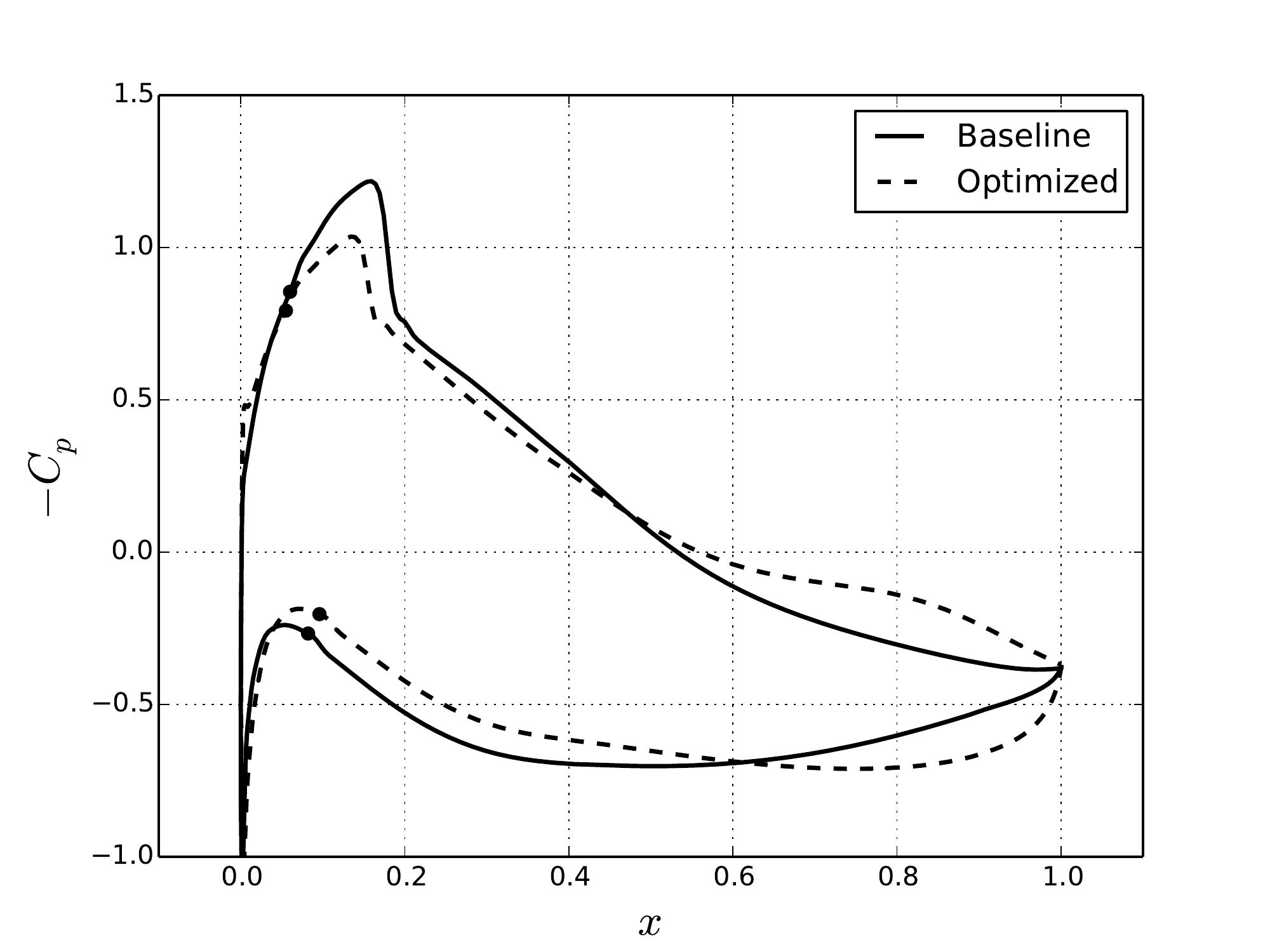}}
\subfigure[$\alpha = 4.5^\circ$]{\includegraphics[width=0.49\textwidth]{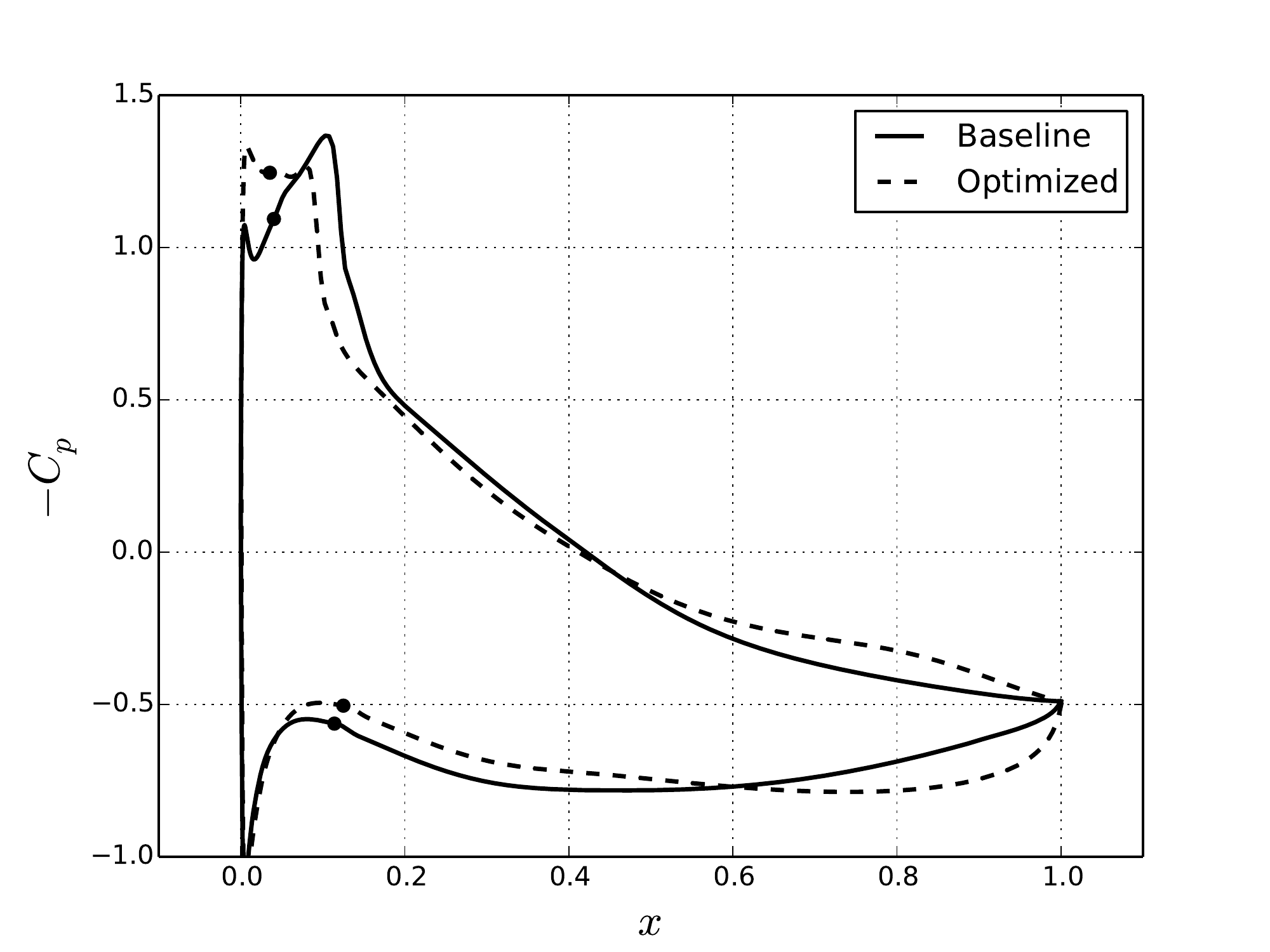}}
\caption{Comparison between the baseline and multi-point deterministic optimal pressure coefficient profiles for the UTRC blade at three different incidence angles. The transition locations are indicated by dots.}
\label{fig:uti_cp_comb_mp}
\end{figure*}

The pressure coefficient profiles for the baseline and multi-point deterministic optimal geometries are compared in Figure \ref{fig:uti_cp_comb_mp}, showing that the optimizer shifts the loading toward the trailing edge which reduces the strength of the suction side shock. At the negative incidence design point, reducing the shock strength reduces the nominal and mean loss on the suction side of the blade. Flow separation occurs on the suction surface for only 5\% of the optimized manufactured blades, whereas suction side flow separation occurs on 35\% of the manufactured baseline blades. On the pressure side of the blade, the design changes have the largest impact at negative incidence. Decreasing the blade stagger reduces the number of blades with pressure side flow separation from 50\% down to 14\%.

The single-point deterministic optimal design eliminated the shock at the nominal incidence, but increased the shock strength at positive incidence relative to the baseline design. The multi-point deterministic optimal does not eliminate the suction side shock at any of the design points, so a switch in the dominant loss mechanism does not occur when manufacturing variations are introduced. Including variations in incidence has a similar effect as including manufacturing variations: the multi-point optimal design avoids regions of the design space where a switch in the dominant loss mechanism occurs. For the multi-point optimization, robustness to flow incidence implies robustness to manufacturing variations, since changes in incidence have a larger effect on the mean performance than manufacturing variability.

\begin{figure*}[htbp]
\centering
\subfigure[Deterministic Optimal Tolerances]{\includegraphics[width=0.49\textwidth]{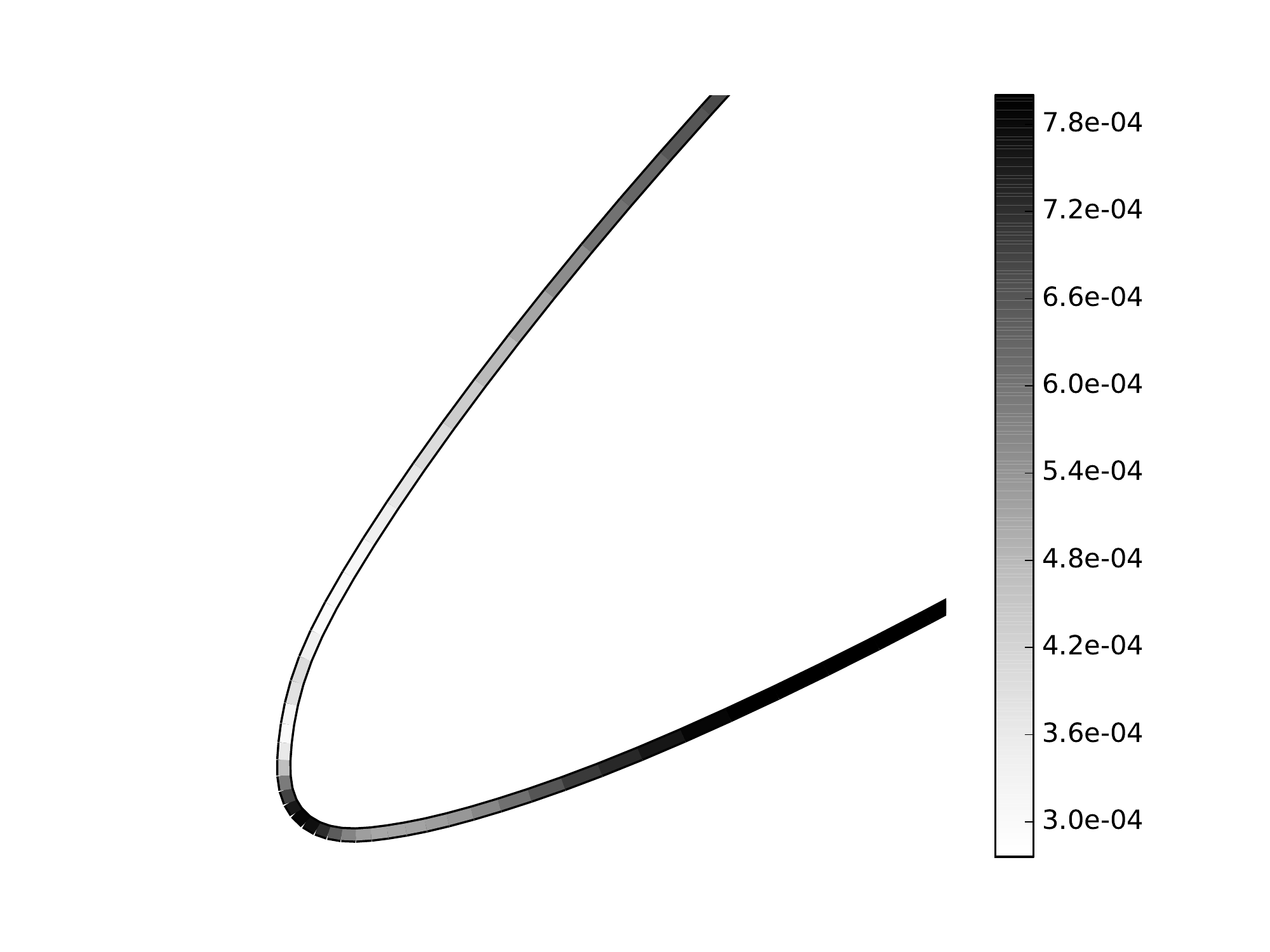}}
\subfigure[Robust Optimal Tolerances]{\includegraphics[width=0.49\textwidth]{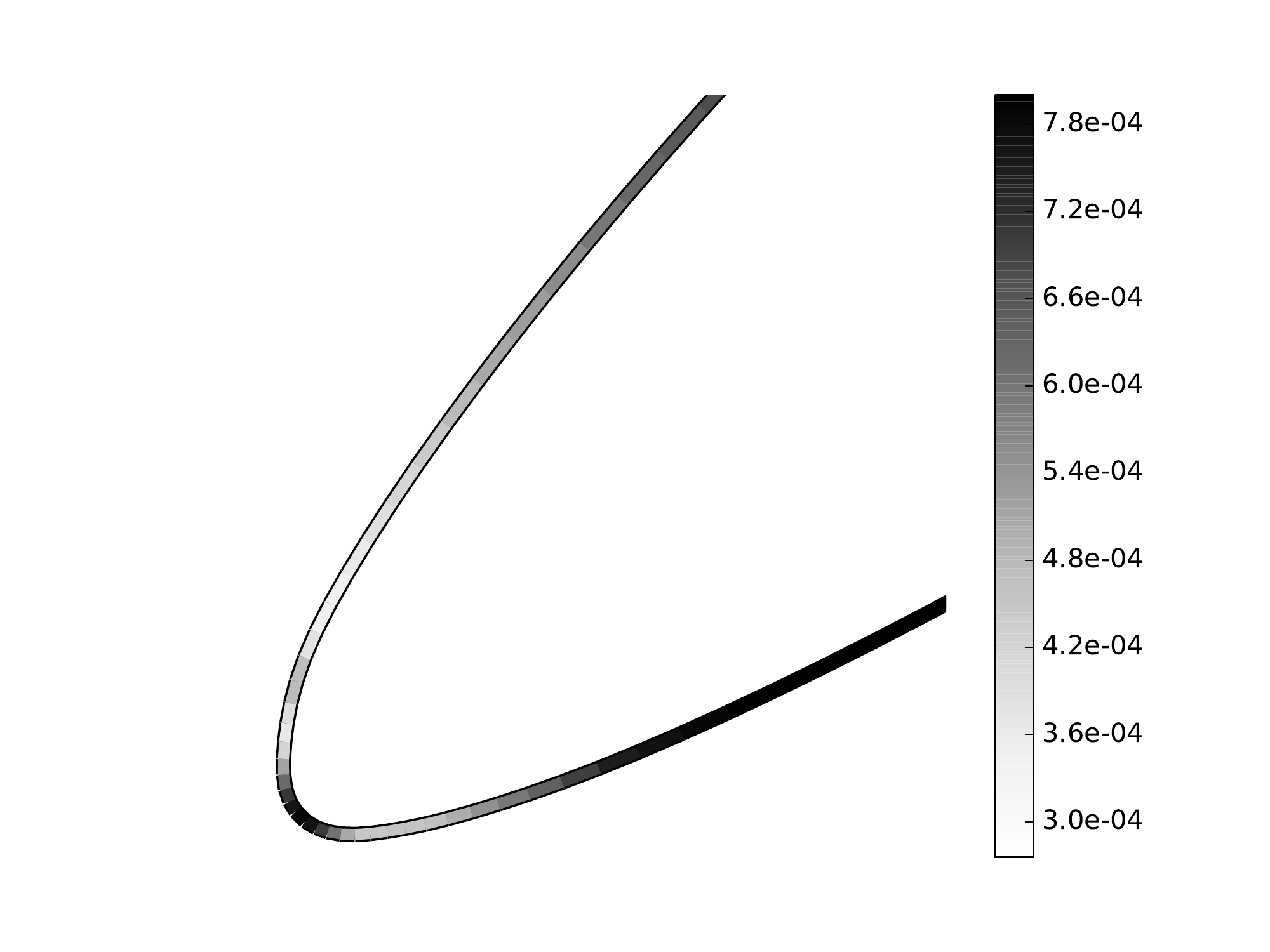}}
\subfigure[Simultaneous Optimal Tolerances]{\includegraphics[width=0.49\textwidth]{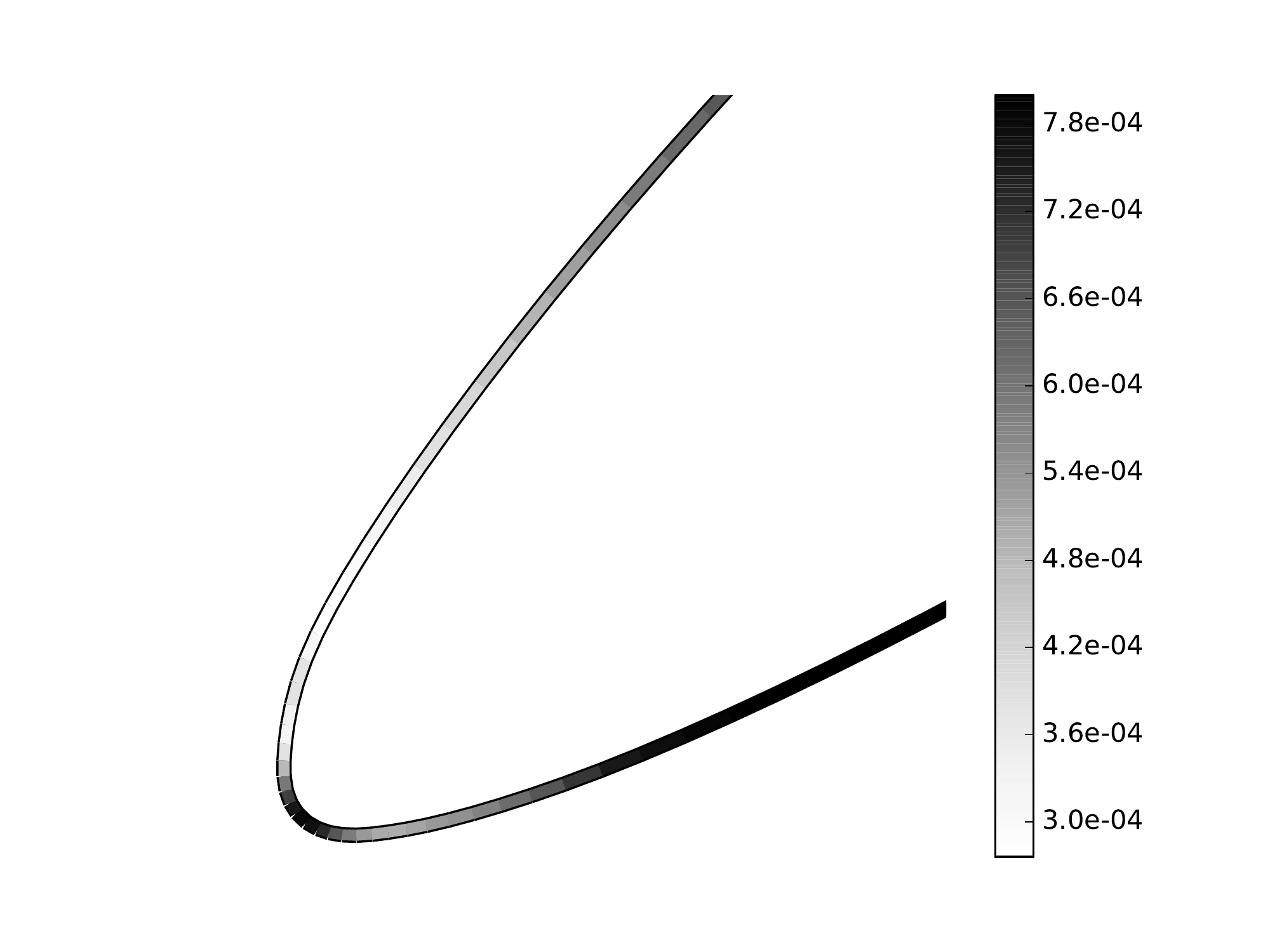}}
\caption{Optimal standard deviation $\sigma(s)/c$ for the multi-point optimized UTRC blades. The lower surface is the pressure side, and the upper surface is the suction side.}
\label{fig:uti_mp_opt_tol}
\end{figure*}

The insensitivity of the optimal multi-point geometry to manufacturing variations implies that sequential tolerance design should give similar results to tolerances produced by the simultaneous approach. The largest error between the optimal tolerances, plotted in Figure \ref{fig:uti_mp_opt_tol} is 8\%, roughly four times smaller than the maximum error between the single-point optimal tolerances. The multi-point optimal tolerances decrease the standard deviation near the leading edge of the blade, resulting in a ``double bow-tie'' tolerance pattern. The greatest decrease is on the suction side of the blade, since most of the mean shift arises due to suction side separation occurring at positive incidence.

\section{Conclusions and Future Work} 

This paper has presented a framework for designing the distribution of manufacturing variability for turbomachinery blades. Applying this framework to a subsonic exit guide vane geometry produced a novel ``double bow-tie'' tolerance scheme that reduced the extent of flow separation away from the design incidence angle. Geometry optimization was incorporated into the tolerance design framework to optimize both the design intent blade geometry and manufacturing tolerances. Single-point deterministic geometry optimization produced a design that was neither robust to manufacturing variations nor to variations in flow incidence. A switch in the dominant loss mechanism degrades the performance of the single-point deterministic optimal designs when manufacturing variations are introduced. Multi-point deterministic optimization produced optimal blade geometries that were robust with respect to flow incidence and manufacturing variations, provided the design variables modify the geometry on length scales larger than the length scales of the manufacturing variations. Such design variables include those that alter the camber and stagger of the blade. Similar behavior was observed when optimizing a subsonic rotor blade \cite{dow_thesis}. When modifying the geometry at smaller length scales, the design procedures recommended by Goodhand \textit{et al.} in \cite{goodhand_2014} should be adopted.

In practice, the proposed framework for tolerance design can be used to select manufacturing processes that achieve the required process spread at every location on the blade, and to construct tolerance bands for rejecting or reworking blades. The proposed approach determines the optimal process spread $\sigma(s)$, which varies continuously around the blade. In reality, only finitely many manufacturing processes are available, so $\sigma(s)$ changes discontinuously at the point where one manufacturing process is replaced by another. Further work is required to determine how to translate the continuous process spread to a set of manufacturing processes in different regions of the blade.

Tolerance bands can also be constructed from the process spread by choosing constant factors $k_i$ and $k_o$ corresponding to inner and outer tolerance bands. Blades whose surface falls within the inner tolerance band $x_d(s) - k_i\sigma(s)\hat{n}(s)$ are rejected, and blades whose surface falls outside the outer tolerance band $x_d(s) + k_o\sigma(s)\hat{n}(s)$ are reworked. Rejecting and reworking blades that fall outside the tolerance bands truncates the distribution of the manufactured blade variability, as opposed to scaling the distribution according to the process spread. Future work should recommend how to select $k_i$ and $k_o$, and compare the mean performance of blades toleranced according to bands constructed from $\sigma(s)$ to the mean performance of blades with scaled process spread.

%The multi-point robust optimal compressor blades are robust to geometric variability introduced by the manufacturing process and to incidence variations. Future work should investigate the impact of variability in the inflow Mach number on the robust optimal design. Variability in the inflow Mach number can be incorporated to the design framework using Monte Carlo to sample from a distribution of Mach numbers as was done to incorporate variations in the geometry, or using quadrature as was done to incorporate variations in the flow incidence. Results presented in Chapter \ref{chap:chap3} for the single-point optimized exit guide vane showed that the formation of shocks can cause a switch in the dominant loss mechanism of manufactured blades, and the same switch could be triggered by variations in the inflow Mach number. Future work should address this possibility and recommend design practices in such cases.

%%%%%%%%%%%%%%%%%%%%%%%%%%%%%%%%%%%%%%%%%%%%%%%%%%%%%%%%%%%%%%%%%%%%%%

%%%%%%%%%%%%%%%%%%%%%%%%%%%%%%%%%%%%%%%%%%%%%%%%%%%%%%%%%%%%%%%%%%%%%%

\begin{acknowledgment}
The authors would like to thank Professor Mark Drela for his suggestions on running MISES. The authors would also like to thank Professor Dave Darmofal, Professor Rob Miller, Professor Ed Greitzer and Professor Karen Willcox for their discussion and feedback on the results. Financial support from Pratt \& Whitney and The Boeing Company is gratefully acknowledged.
\end{acknowledgment}

%%%%%%%%%%%%%%%%%%%%%%%%%%%%%%%%%%%%%%%%%%%%%%%%%%%%%%%%%%%%%%%%%%%%%%
% The bibliography is stored in an external database file
% in the BibTeX format (file_name.bib).  The bibliography is
% created by the following command and it will appear in this
% position in the document. You may, of course, create your
% own bibliography by using thebibliography environment as in
%
% \begin{thebibliography}{12}
% ...
% \bibitem{itemreference} D. E. Knudsen.
% {\em 1966 World Bnus Almanac.}
% {Permafrost Press, Novosibirsk.}
% ...
% \end{thebibliography}

% Here's where you specify the bibliography style file.
% The full file name for the bibliography style file 
% used for an ASME paper is asmems4.bst.
\bibliographystyle{asmems4}

% Here's where you specify the bibliography database file.
% The full file name of the bibliography database for this
% article is asme2e.bib. The name for your database is up
% to you.
\bibliography{asme2ej}

\begin{thebibliography}{10}

\bibitem{garzon_thesis}
Garzon, V.~E., 2003.
\newblock ``Probabilistic aerothermal design of compressor airfoils''.
\newblock {PhD} dissertation, Massachusetts Institute of Technology, Department
  of Aeronautics and Astronautics.

\bibitem{kumar_2006_2}
Kumar, A., Keane, A., Nair, P., and Shahpar, S., 2005.
\newblock ``Robust design of compressor fan blades against erosion''.
\newblock {\em Journal of Mechanical Design, {\bf 128}}, pp.~864--873.

\bibitem{bestle_2010}
Bestle, D., and Flassig, P., 2010.
\newblock ``Optimal aerodynamic compressor blade design considering
  manufacturing noise''.
\newblock In Proceedings of the 8th {ASMO UK/ISSMO} Conference.

\bibitem{goodhand_2014}
Goodhand, M., Miller, R., and Lung, H., 2014.
\newblock ``The impact of geometric variation on compressor 2d incidence
  range''.
\newblock {\em Journal of Turbomachinery}.

\bibitem{beachkofski_2004}
Beachkofski, B., and Grandhi, R., 2004.
\newblock ``Probabilistic system reliability for a turbine engine airfoil''.
\newblock In Proceedings of ASME Turbo Expo 2004.

\bibitem{duffner_thesis}
Duffner, J., 2008.
\newblock ``The effects of manufacturing variability on turbine vane
  performance''.
\newblock {Master's} dissertation, Massachusetts Institute of Technology,
  Department of Aeronautics and Astronautics.

\bibitem{sinha_2008}
Sinha, A., Hall, B., Cassenti, B., and Hilbert, G., 2008.
\newblock ``Vibratory parameters of blades from coordinate measurement machine
  data''.
\newblock {\em Journal of Turbomachinery, {\bf 130}}.

\bibitem{lange_2012}
Lange, A., Voigt, M., Vogeler, K., and Johann, E., 2012.
\newblock ``Principal component analysis of {3D} scanned compressor blades for
  probabilistic {CFD} simulation''.
\newblock In Proceedings of the 53rd AIAA/ASME/ASCE/AHS/ASC Structures,
  Structural Dynamics and Materials Conference.

\bibitem{loeve}
Lo{\`e}ve, M., 1963.
\newblock ``Probability theory''.
\newblock {\em Graduate Texts in Mathematics, {\bf 45}}.

\bibitem{lemaitre}
Le~M{\^a}itre, O., and Knio, O., 2010.
\newblock {\em Spectral Methods for Uncertainty Quantification - With
  Applications to Computational Fluid Dynamics}, 1st~ed.
\newblock Springer Verlag, New York, ch.~2, pp.~17--44.

\bibitem{kane_1986}
Kane, V., 1986.
\newblock ``Probabilistic optimal design in the presence of random fields''.
\newblock {\em Journal of Quality Technology, {\bf 18}}, pp.~41--52.

\bibitem{rubinstein_shapiro}
Rubinstein, R., and Shapiro, A., 1993.
\newblock {\em Discrete event systems: Sensitivity analysis and stochastic
  optimization by the score function method}, Vol.~346.
\newblock Wiley, New York.

\bibitem{powell_1978}
Powell, M., 1978.
\newblock ``A fast algorithm for nonlinearly constrained optimization
  calculations''.
\newblock {\em Lecture Notes in Mathematics, {\bf 630}}, pp.~144--157.

\bibitem{stephens_hobbs_1979}
Stephens, H.~E., and Hobbs, D.~E., 1979.
\newblock Design and performance evaluation of supercritical airfoils for axial
  flow compressors.
\newblock Tech. Rep. PWA-FR-11455, Pratt and Whitney Government Products
  Division.

\bibitem{camp_1995}
Camp, T., and Shin, H., 1995.
\newblock ``Turbulence intensity and length scale measurements in multistage
  compressors''.
\newblock {\em Journal of Turbomachinery, {\bf 117}}, pp.~38--46.

\bibitem{mises_manual}
Drela, M., and Youngren, H., 2008.
\newblock {\em A User's Guide to {MISES} 2.63}.
\newblock MIT Aerospace Computational Design Laboratory, 70 Vassar St,
  Cambridge MA 02139, February.

\bibitem{garzon_2003}
Garzon, V.~E., and Darmofal, D., 2003.
\newblock ``Impact of geometric variability on axial compressor performance''.
\newblock {\em Journal of Turbomachinery, {\bf 125}}(4), pp.~692--703.

\bibitem{dow_thesis}
Dow, E.~A., 2014.
\newblock ``Robust design and tolerancing of compressor blades''.
\newblock {PhD} dissertation, Department of Aeronautics and Astronautics.

\end{thebibliography}

%%%%%%%%%%%%%%%%%%%%%%%%%%%%%%%%%%%%%%%%%%%%%%%%%%%%%%%%%%%%%%%%%%%%%%
%\appendix       %%% starting appendix
%\section*{Appendix A: Head of First Appendix}
%Avoid Appendices if possible.

\end{document}